\documentclass[12pt]{article}
\usepackage{amssymb}
\usepackage{amsmath}
\usepackage{indentfirst}
\usepackage{graphicx}
\usepackage{mathrsfs}

\def\squarebox#1{\hbox to #1{\hfill\vbox to #1{\vfill}}}
\newcommand{\qedbox}{\vbox{\hrule\hbox{\vrule\squarebox{.667em}\vrule}\hrule}}
\newcommand{\qed}{\nopagebreak\mbox{}\hfill\qedbox\bigskip\par}

\def\Proof{\removelastskip\vskip\baselineskip\noindent{\bf
Proof\quad}\ignorespaces}

\begin{document}
\title{Combinatorial properties of the numbers of tableaux of bounded height}
\date{}\author{Marilena Barnabei, Flavio Bonetti, and Matteo
Silimbani
\thanks{ Department of Mathematics - University of Bologna}}
\maketitle

\begin{abstract}
\noindent We introduce an infinite family of lower triangular
matrices $\Gamma^{(s)}$, where $\gamma_{n,i}^s$ counts the
standard Young tableaux on $n$ cells and with at most $s$ columns
on a suitable subset of shapes. We show that the entries of these
matrices satisfy a three-term row recurrence and we deduce
recursive and asymptotic properties for the total number
$\tau_s(n)$ of tableaux on $n$ cells and with at most $s$ columns.
\end{abstract}

\section{Introduction}

\noindent The first simple expressions of the number of standard
Young tableaux of given shape were given by Frobenius and Young
(see \cite{fy} and \cite{yf}) and by Frame-Robinson-Thrall
\cite{frt}. More recently, the number of standard Young tableaux
has been studied according to the height of their shape. Regev
\cite{reg} gave asymptotic values for the numbers $\tau_s(n)$ of
standard Young tableaux whose shape consists of $n$ cells and at
most $s$ columns and Stanley (see \cite{sta}) discussed the
algebraic or differentiably finite nature of the corresponding
generating functions. More recently, many authors entered this
vein. For example, D.Gouyou-Beauchamps gave both exact forms in
\cite{gb} and recurrence formulas (see \cite{gbphd}) for the
numbers $\tau_s(n)$ when $s\leq 5$ by combinatorial tools, while
F.Bergeron and F.Gascon \cite{bg} found some recurrence formulas
for the numbers $\tau_s(n)$ by analytic argumentations.

\noindent In this paper we present a different approach to the
subject. Our starting point is the remark that the elements of the
Ballot matrix introduced by M.Aigner in \cite{aig} correspond
bijectively to the integers counting standard Young tableaux of a
given shape with at most $2$ columns. Firstly, we arrange the
entries of the Ballot Matrix in a new lower triangular matrix $A$
in such a way that the entries of the $n$-th row count standard
Young tableaux with precisely $n$ cells. The integer $\tau_2(n)$
can be therefore recovered by summing the elements of the $n$-th
row of $A$.

\noindent In the following sections, we extend the results to the
general case of standard Young tableaux with at most $s$ columns.
For every fixed $s\in\mathbb{N}$ we define an infinite matrix
$\Gamma^{(s)}$ whose $(n,i)$-th entry is the total number of
standard Young tableaux with at most $s$ columns and such that the
difference between the length of the second and of the third
column is $i$. The total number $\tau_s(n)$ of tableaux on $n$
cells with at most $s$ columns is the $n$-th row sum of the matrix
$\Gamma^{(s)}$. The matrix $\Gamma^{(s)}$ presents a three-term
row-recurrence property, which yields a recurrence law satisfied
by the integers $\tau_s(n)$. This recurrence allows to get a
combinatorial interpretation of the asymptotic behaviour of the
ratio
$$\frac{\tau_s(n)}{\tau_s(n-1)}.$$
In Section \ref{ter} we treat more extensively the case $s=3$, in
order to exhibit recursive and asymptotic properties for the
integers $\tau_3(n)=M_n$, where $M_n$ denotes the $n$-th Motzkin
number.

\section{The two-column case}
\label{sec}

\noindent In the following, we will denote the shape of a standard
Young tableau $T$ as a list containing the length of the first,
second, $\ldots$, last column of $T$.

\noindent We first consider the standard Young tableaux whose
shape consists of at most two columns. We define an infinite,
lower triangular matrix $A=(\alpha_{n,i})$, where $\alpha_{n,i}$
is the number of standard Young tableaux whose shape corresponds
to the list $(n-i,i)$, with $i\in\mathbb{N}$.

\noindent The matrix $A$ can be recursively constructed as
follows:

\newtheorem{recdef}{Proposition}

\begin{recdef}
The matrix $A$ is defined by the following initial conditions:
\begin{itemize}\label{inverso}
\item[i)] $\alpha_{i,0}=1$ for every $i\geq 0$;
\item[ii)] $\alpha_{1,i}=0$ for every $i>0$;
\item[iii)] $\alpha_{n,i}=0$ for every $i>\lfloor\frac{n}{2}\rfloor$
\end{itemize}
and by the recurrence: \begin{equation}
\alpha_{n,i}=\alpha_{n-1,i}+\alpha_{n-1,i-1} \textrm{\ \  for\ \
} n\geq 1, 1\leq i\leq
\lfloor\frac{n}{2}\rfloor.\label{iv}\end{equation}

\end{recdef}
\Proof The initial conditions follow immediately by the definition
of the matrix $A$. In order to prove the recurrence, remark that a
Young tableau with $n$ cells can be obtained from a tableau with
$n-1$ cells by adding a new box containing the symbol $n$ in a
corner position. In particular, in a Young tableau of shape
$(n-i,i)$, $i>0$, the symbol $n$ can be placed in at most two
corner cells. Hence, the entry $\alpha_{n,i}$ is the sum of the
two integers $\alpha_{n-1,i-1}$ and $\alpha_{n-1,i}$, that count
the Young tableaux of appropriate shape.\qed

\noindent The recurrence formula (\ref{iv}) has as an immediate
fallout the following result:

\newtheorem{recdue}[recdef]{Corollary}

\begin{recdue}
The entries of the matrix $A$ satisfy the following columnwise
recurrence:

$$\alpha_{n,i}=\sum_{h=2i-1}^{n-1} \alpha_{h,i-1}$$ for every
$n,i>0$.

\end{recdue}

\noindent We remark that the elements of the matrix $A$ correspond
bijectively to the entries of the Ballot Matrix defined by Aigner
in \cite{aig}. More precisely, the Ballot Matrix
$\widetilde{A}=(\tilde{\alpha}_{i,j})$ can be obtained rearranging
the entries of the matrix $A$ as follows:
$$\tilde{\alpha}_{j,k}=\alpha_{2j-k,j-k}.$$
As a consequence, for every $n$, the number of standard Young
tableaux with two columns of the same length $n$ is the $n$-th
Catalan number $C_n$:
\begin{equation}\alpha_{2n,n}=C_n\label{cut}\end{equation}.
\newline

\noindent Denote by $\tau_2(n)$ the number of Young tableaux with
$n$ cells and at most two columns. Obviously, we have:
$$\tau_2(n)=\sum_{i=0}^{\left\lfloor\frac{n}{2}\right\rfloor}a_{n,i}.$$
In the next theorem we deduce a recurrence formula for the
sequence $\tau_2(n)_{n\geq 1}$:

\newtheorem{catrec}[recdef]{Theorem}

\begin{catrec}
The integers $\tau_2(n)$, $n\geq 1$, satisfy the recurrence
\begin{equation}\tau_2(n)=2\tau_2(n-1)-E(n-1)\cdot C_{
\frac{n-1}{2}},\label{peruno}\end{equation} where $C_i$ denotes
the $i$-th Catalan number and
$$E(n)=\Big\{\begin{array}{ccccc}
1&&&\textrm{if n is even}\\
0&&&\textrm{otherwise.}
\end{array}$$
\end{catrec}

\Proof By Proposition \ref{inverso}, we have:
$$\tau_2(n)=\sum_{i=0}^{\left\lfloor\frac{n}{2}\right\rfloor}a_{n,i}=
\sum_{i=0}^{\left\lfloor\frac{n}{2}\right\rfloor}a_{n-1,i-1}+\sum_{i=0}^{\left\lfloor\frac{n}{2}\right\rfloor}a_{n-1,i}
=$$$$=2\tau_2(n)-E(n-1)\cdot\alpha_{n-1,\frac{n-1}{2}}
=2\tau_2(n)-E(n-1)\cdot C_{ \frac{n-1}{2}},$$ as desired. \qed

\begin{figure}[ht]
\begin{center}
\includegraphics[bb=125 632 266 684,width=.5\textwidth]{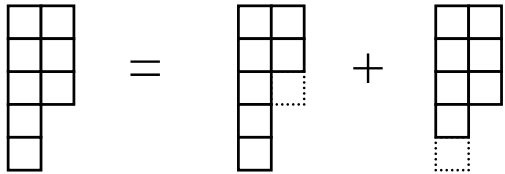} \caption{The
recurrence formula for the entry
$\alpha_{8,3}$.}\label{disegnonoa}
\end{center}
\end{figure}

\noindent The preceding result yields an alternative proof of the
following theorem, originally proved by Regev in \cite{reg}:

\newtheorem{princdue}[recdef]{Theorem}

\begin{princdue}
The number of standard Young tableaux with exactly $n$ cells and
at most two columns is
$$\tau_2(n)={n \choose \lfloor
\frac{n}{2}\rfloor}.$$
\end{princdue}

\Proof It is well known that the central binomial coefficients ${n
\choose \lfloor \frac{n}{2}\rfloor}$ satisfy the recurrence
(\ref{peruno}). Remarking that $\tau_2(0)=1={0 \choose 0}$, we
have the assertion. \qed

\noindent Moreover, identity (\ref{peruno}) yields the following
asymptotic property of the integers $\tau_2(n)$:

\newtheorem{nonlo}[recdef]{Proposition}
\begin{nonlo}\label{becchi}
We have:
$$\lim_{n\to\infty}\frac{\tau_2(n)}{\tau_2(n-1)}=2$$
\end{nonlo}

\Proof By (\ref{peruno}) we get:
$$\frac{\tau_2(n)}{\tau_2(n-1)}=\frac{2\tau_2(n-1)-E(n+1)\cdot C_{\frac{n-1}{2}}}{\tau_2(n-1)}.$$
This implies that the assertion is proved as soon as we show that
$$\lim_{n\to\infty}\frac{E(n+1)\cdot C_{\frac{n-1}{2}}}{\tau_2(n-1)}=0.$$ If $n$ is even, this
identity trivially holds. If $n$ is odd, the described correction
term is negligeable with respect to $\tau_2(n-1)$ as $n$ goes to
infinity. In fact, we get:
$$\lim_{n\to\infty}\frac{ \frac{1}{n}\Big({n-1 \atop
\frac{n-1}{2}}\Big)}{\Big({n-1 \atop
\frac{n-1}{2}}\Big)}=\lim_{n\to\infty}\frac{1}{n}=0,$$ as
required. \qed

\section{The three-column case}\label{ter}

\noindent We now extend the results of Section \ref{sec} to the
case of standard Young tableaux with at most three columns. Define
a matrix $B=(\beta_{n,i})$ such that the $(n,i)$-th entry of $B$
is the total number of standard Young tableaux of shape
$(n-i-2k,i+k,k)$, for every possible value of $k$. In other terms,
the integer $\beta_{n,i}$ is the cardinality of the set
$Y^3_{n,i}$ containing all standard tableaux on $n$ cells with at
most $3$ columns such that the difference between the second and
third column is $i$.

\noindent The entries of the matrix $B$ satisfy a three-term
recurrence:

\newtheorem{dsw}[recdef]{Proposition}

\begin{dsw}\label{chegli}
The integers $\beta_{n,0}$ satisfy the identity:
\begin{equation}\beta_{n,0}=\beta_{n-1,0}+\beta_{n-1,1}.\label{izero}\end{equation} Moreover, for every $1\leq
i\leq\lfloor\frac{n}{2}\rfloor$, we have:
\begin{equation}\beta_{n,i}=\beta_{n-1,i-1}+\beta_{n-1,i}+\beta_{n-1,i+1}-r_{n,i},\label{igrande}\end{equation}
where $$r_{n,i}=\left\{\begin{array}{ll}  \textrm{number of
tableaux of shape }
(\frac{n+i-2}{3},\frac{n+i-2}{3},\frac{n-2i+1}{3})&\textrm{ if }
n-2i\equiv_3 2\\ 0 & otherwise.\end{array}\right.$$
\end{dsw}

\Proof Let $T$ be a standard Young tableau $T$ with $3$ columns of
shape $(n-i-2k,i+k,k)$, for some $k$, and let $T'$ be the tableau
obtained from $T$ by removing the cell containing the symbol $n$.
Obviously, the tableau $T'$ has either shape $(n-i-2k-1,i+k,k)$ or
$(n-i-2k,i+k-1,k)$ or $(n-i-2k,i+k,k-1)$. This implies that the
correspondence $T\mapsto T'$ maps the set $Y^3_{n,i}$ into the
union $Y^3_{n-1,i-1}\cup Y^3_{n-1,i}\cup Y^3_{n-1,i+1}$. However,
such map is not in general surjective. In fact, it is easy to
check that the difference between the set $Y^3_{n-1,i-1}\cup
Y^3_{n-1,i}\cup Y^3_{n-1,i+1}$ and the image of the map consists
exactly of the tableaux of shape
$(\frac{n+i-2}{3},\frac{n+i-2}{3},\frac{n-2i+1}{3})$. This proves
(\ref{igrande}). The first recurrence can be proved by the same
arguments, observing that, in this case, the map defined above is
surjective.\qed

\noindent The preceding result yields a recurrence formula
satisfied by the total number $\tau_3(n)$ of standard Young
tableaux with $n$ cells and at most $3$ columns:

\newtheorem{casino}[recdef]{Theorem}

\begin{casino}\label{ultcor}
The numbers $\tau_3(n)$, $n\geq 3$, satisfy the recurrence
\begin{equation}\tau_3(n)=3\tau_3(n-1)-R(n)-E(n-1)\cdot C_{ \frac{n-1}{2}}-\beta_{n-1,0},\label{nsbpn}\end{equation}where
$$R(n)=\sum_{i=0}^{\left\lfloor\frac{n}{2}\right\rfloor} r_{n,i}.$$
\end{casino}

\Proof  By definition of $\beta_{n,i}$ we have:
$$\tau_3(n)=
\sum_{i=0}^{\left\lfloor\frac{n}{2}\right\rfloor}\beta_{n,i}.$$
Hence,
\begin{equation}
\tau_3(n)=\sum_{i=0}^{\left\lfloor\frac{n}{2}\right\rfloor}\beta_{n-1,i-1}+\sum_{i=0}^{\left\lfloor\frac{n}{2}\right\rfloor}\beta_{n-1,i}
+\sum_{i=0}^{\left\lfloor\frac{n}{2}\right\rfloor}\beta_{n-1,i+1}-
\sum_{i=0}^{\left\lfloor\frac{n}{2}\right\rfloor}r_{n-1,i-1}.\label{rihs}\end{equation}
Remark that
$$\tau_3(n-1)=\sum_{i=0}^{\left\lfloor\frac{n-1}{2}\right\rfloor}\beta_{n-1,i}.$$
The proof is an immediate consequence of the following
considerations:
\begin{itemize}
\item the first summand in the right hand side of Identity
(\ref{rihs}) is equal to $\tau_3(n-1)$ if $n$ is even. In case of
$n$ odd, instead, this term exceeds $\tau_3(n-1)$ by the Catalan
number $C_{\frac{n-1}{2}}$,
\item the second summand in the right hand side of Identity
(\ref{rihs}) is exactly $\tau_3(n-1)$,
\item the third summand in the right hand side of Identity
(\ref{rihs}) exceeds the integer $\tau_3(n-1)$ exactly by
$\beta_{n-1,0}$,
\item the fourth summand in the right hand side of Identity
(\ref{rihs}) is obviously equal to $R(n-1)$.
\end{itemize}
\qed

\noindent Remark that an explicit evaluation of the integers
$Q(n)$ and $R(n)$ can be obtained via the well known Hook Length
Formula.\newline

\noindent Previous considerations lead to the following result
concerning the asymptotic behaviour of the sequence $\tau_3(n)$:

\newtheorem{aigcomb}[recdef]{Proposition}

\begin{aigcomb}
\label{aneso}For every integer $n\geq 1$, we have:
\begin{equation}\frac{\tau_3(n)}{\tau_3(n-1)}<
3.\label{maserv}\end{equation}
 Moreover\textrm{,}
\begin{equation}\lim_{n\to\infty}\frac{\tau_3(n)}{\tau_3(n-1)}=3.\label{siche}\end{equation}
\end{aigcomb}

\Proof Inequality (\ref{maserv}) follows directly by previous
considerations. In order to prove (\ref{siche}) we remark that, by
Theorem (\ref{ultcor}), the ratio $\frac{\tau_(n)}{\tau_(n-1)}$
can be written as:
$$\frac{\tau_3(n)}{\tau_3(n-1)}=3-(U_1(n-1)+U_2(n-1)+U_3(n-1)),$$ where
$$U_1(n-1)=\frac{E(n+1)\cdot C_{\frac{n-1}{2}}}{\tau_3(n-1)}$$
$$U_2(n-1)=\frac{\beta_{n-1,0}}{\tau_3(n-1)}$$
$$U_3(n-1)=\frac{R(n-1)}{\tau_3(n-1)}.$$
The statement will be proved as soon as we show that each summand
$U_i(n-1)$ goes to $0$ as $n$ goes to infinity. The numerator of
each $U_i(n-1)$ is a sum of a suitable finite number
$f_{\lambda}$. For every $f_{\lambda}$ appearing in one of these
expression, we can single out an appropriate subset $S_{\lambda}$
of summands whose cardinality grows as $n$ goes to infinity, such
that:
$$\lim_{n\to\infty}\frac{f_{\lambda}}{\tau_{3}(n-1)}\leq\lim_{n\to\infty}\frac{f_{\lambda}}{\sum_{\mu\in S_{\lambda}}f_{\mu}}= 0.$$
For instance, consider the rectangular shape $\rho$ over $n-1$
blocks, namely, the shape corresponding to the values $i=0$,
$k=\frac{n-1}{3}$. Obviously, such a shape exists only if
$n-1\equiv 0$ (mod 3). In this case, the family $S_{\rho}$
consists of the shapes over $n-1$ cells and corresponding to the
values $i=h$, $k=k=\frac{n-1}{3}-h$, where $h$ ranges from $1$ to
$k=\frac{n-1}{6}$. By Hook Length Formula, we get:
$$\lim_{n\to\infty}\frac{f_{\rho}}{\tau_3(n-1)}\leq\lim_{n\to\infty}\frac{f_{\rho}}{\sum_{\mu\in
S_{\rho}}f_{\mu}}=\lim_{n\to\infty}\frac{f_{\rho}}{\sum_{h=1}^{\left\lfloor\frac{n-1}{6}\right\rfloor}(h+1)^3f_{\rho}}
\leq\lim_{n\to\infty}\frac{1}{8^{\left\lfloor\frac{n-1}{6}\right\rfloor}}=
0.$$  The sets $S_{\nu}$ corresponding to the other summand of
each of the numerators of the $U_i(n-1)$ can be described
analogously.
 \qed

\noindent As proved in \cite{reg}, $\tau_3(n)=M_n$, where $M_n$ is
the $n$-th Motzkin number. Hence, these last arguments yield a
combinatorial proof of Propositions $4$ and $5$ in \cite{aig2}.

\section{The general case}
\noindent In this section we extend the previous argumentations to
the general case of standard Young tableaux with at most $s$
columns. Define a matrix $\Gamma^{(s)}=(\gamma^{(s)}_{n,i})$ such
that the $(n,i)$-th entry is the total number of standard Young
tableaux with at most $s$ columns and such that the difference
between the length of the second and of the third column is $i$.

\noindent Also in this case, when $s\geq 4$ and $n\geq s$, the
entries of this matrix satisfy a three-term recurrence. In fact,
the following proposition can be proved by the same arguments used
to prove Proposition \ref{chegli}:

\newtheorem{casoqua}[recdef]{Proposition}

\begin{casoqua}
The integers $\gamma^{(s)}_{n,0}$ satisfy the identity:
\begin{equation}\gamma^{(s)}_{n,0}=(s-2)\gamma^{(s)}_{n-1,0}+\gamma^{(s)}_{n-1,1}-\sum_{j=3}^{s-1}r^{(s)}_j(n-1,0).\label{i_zero}\end{equation}
Moreover, for every $1\leq i\leq\lfloor\frac{n}{2}\rfloor$, we
have:
$$\gamma^{(s)}_{n,i}=\gamma^{(s)}_{n-1,i-1}+(s-2)\gamma^{(s)}_{n-1,i}+\gamma^{(s)}_{n-1,i+1}-r^{(s)}_1(n-1,i-1)$$
\begin{equation}-\sum_{j=3}^{s-1}r^{(s)}_j(n-1,i)
,\label{i_grande}\end{equation} where $r^{(s)}_j(n-1,i)$ is the
number of standard Young tableaux with at most $s$ columns such
that the difference between the length of the second and of the
third column is $i$ and the $j$-th and the $(j+1)$-th columns have
the same length.
\end{casoqua}\qed

\noindent Similarly, exploiting the same arguments as in
Proposition \ref{aneso}, we get:

\newtheorem{sfm}[recdef]{Proposition}

\begin{sfm}
\noindent The numbers $\tau_s(n)$ satisfy the recurrence:
\begin{equation}\tau_s(n)=s\tau_s(n-1)-E(n-1)C_{\frac{n-1}{2}}-\gamma^{(s)}_{n-1,0}-
R^{(s)}(n-1),\label{gramat}\end{equation} where $R^{(s)}(n-1)$ is
the sum of all the correction terms appearing in Formula
\emph{(\ref{i_grande})}.
\end{sfm}\qed

\noindent Now we can easily deduce from these identities some
information concerning the asymptotic behaviour of the sequence
$\tau_s(n)$:

\newtheorem{flavio}[recdef]{Proposition}

\begin{flavio}
The sequence $\tau_s(n)$ satisfies the following properties:
$$\frac{\tau_s(n)}{\tau_s(n-1)}<s$$
$$\lim_{n\to\infty}\frac{\tau_s(n)}{\tau_s(n-1)}=s.$$
\end{flavio}\qed

\section*{Acknowledgments}

\noindent We thank Elisa Pergola and Renzo Pinzani for many useful
conversations.


\begin{thebibliography}{99}

\bibitem{aig} M.Aigner, Catalan and other numbers: a recurrente theme, in {\em Algebraic Combinatorics and Theoretical Computer Science}, H.Crapo and
D.Senato eds., Springer-Verlag, (2001), 347--390.

\bibitem{aig2} M.Aigner, Motzkin numbers,
 {\em Europ. J. Combinatorics} {\bf 19}, (1998), 663--675.

\bibitem{bfk} F.Bergeron, L.Favreau, D.Krob, Conjectures on the enumeration of tableaux of bounded height,
 {\em Discrete Math.} {\bf 139}, (1995), 463--468.

\bibitem{bg} F.Bergeron, F.Gascon, Counting Young tableaux of bounded height,
 \emph{J. Integer Seq.} \textbf{3}, No.1, Art. 00.1.7, 7 p., electronic only (2000)

\bibitem{frt} J.S.Frame, G. De B.Robinson and R.Thrall, The hook graphs of the symmetric group,
\emph{Can. J. Math.}, {\bf 6}, (1954), 316--324.

\bibitem{fy} G.Frobenius, Uber die charaktere de symmetrischen gruppen,
\emph{Preuss. Akad. Wiss. Sitz.}, (1900), 516--534.

\bibitem{g} I.Gessel, Symmetric functions and P-recursiveness,
\emph{Jour. of Comb. Th.}, Series A, {\bf 53}, (1990), 257--285.

\bibitem{gb} D.Gouyou-Beauchamps, Standard Young tableaux of height 4 and 5,
 {\em Europ. J. Combinatorics.} {\bf 10}, (1989), 69--82.

\bibitem{gbphd} D.Gouyou-Beauchamps, Codages par des mots et des chemins: probl\`{e}mes combinatoires et
algorithmiques, Ph.D. Thesis, University of Bordeaux I, 1985.

\bibitem{knu}
D.E.Knuth, ``The art of computer programming: sorting and
searching'', Vol. 3, Addison-Wesley (1998).

\bibitem{nak} T.Nakayama, On some modular properties of irreducible representations of a symmetric
group, I, II, \emph{Jap. J. Math.}, {\bf 17}, (1940), 411--423.

\bibitem{reg} A.Regev, Asymptotic values for degrees associated
with strips of Young tableau, {\em Adv. in Math.}, {\bf 41},
(1981), 115--136.

\bibitem{gb} N.J.A. Sloane, Encyclopedia of integer sequences,
http://www.research.att.com/~njas/sequences/

\bibitem{sta} R.P.Stanley, Differentiability finite power series,
\emph{Europ. J. Comb.}, {\bf 1}, (1980), 175--188.

\bibitem{yf} A.Young, On quantitative substitutional analysis II,
\emph{Proc. Lond. Math. Soc.}, {\bf 34}, (1902), 361--397.




\end{thebibliography}
\end{document}